\documentclass{gtpart}
\usepackage{geometry}               % See geometry.pdf to learn the layout options. There are lots.
\usepackage{graphicx}
\usepackage{nicefrac}
\usepackage{epstopdf}
\usepackage{amsfonts}
\usepackage{amsxtra}     % Use various AMS packages
\usepackage{epsfig}
\usepackage{verbatim}

\title{On volumes of hyperbolic orbifolds}
\author{Ilesanmi Adeboye}
\givenname{Ilesanmi}
\surname{Adeboye}
\address{Department of Mathematics and Computer Science\\
         Wesleyan University\\
        Middletown, CT 06459}
\email{iadeboye@wesleyan.edu}

\author{Guofang Wei}
\givenname{Guofang}
\surname{Wei}
\address{Department of Mathematics\\
         University of California\\
         Santa Barbara, CA 93106}
\email{wei@math.ucsb.edu}


\arxivreference{0911.4712}
\arxivpassword{33rzh}

\DeclareMathOperator{\Vol}{Vol}
\DeclareMathOperator{\Ker}{Ker}
\DeclareMathOperator{\ad}{ad}
\DeclareMathOperator{\spn}{span}
\DeclareMathOperator{\tr}{trace}

\DeclareMathOperator{\Out}{Out}

\DeclareMathOperator{\expo}{exp}

\providecommand{\norm}[1]{\lVert#1\rVert}

 \newcommand{\be}{\begin{equation}}
       \newcommand{\ee}{\end{equation}}
       \newcommand{\ba}{\begin{eqnarray}}
        \newcommand{\ea}{\end{eqnarray}}
 \newcommand{\ban}{\begin{eqnarray*}}
 \newcommand{\ean}{\end{eqnarray*}}

 \newcommand{\lp}{\langle}
 \newcommand{\rp}{\rangle}
 \newcommand{\ra}{\rightarrow}

 \renewcommand{\theequation}{\arabic{section}.\arabic{equation}}
 \newcommand{\sect}[1]{\section{#1} \setcounter{equation}{0}}

 \newcommand{\vol}{\mathrm{Vol}}

\newtheorem{theorem}{Theorem}[section]
\newtheorem{prop}[theorem]{Proposition}
\newtheorem{corollary}[theorem]{Corollary}
\newtheorem{lemma}[theorem]{Lemma}

\newtheorem{definition}[theorem]{Definition}

\newtheorem{remark}[theorem]{Remark}

\date{}
\begin{document}
\begin{abstract} In this paper we use H. C. Wang's bound on the radius of a ball embedded in the fundamental domain of a lattice of a semisimple Lie group to construct an explicit lower bound for the volume of a hyperbolic $n$--orbifold.
\end{abstract}
\bibliographystyle{plain}
\maketitle
\section{Introduction} 
Let $\mathbb{H}^n$ denote hyperbolic $n$--space; the unique simply connected $n$--dimensional Riemannian manifold of constant sectional curvature $-1$.  A \textit{hyperbolic $n$--orbifold} $Q$ is a quotient $\mathbb{H}^n/\Gamma$,  where $\Gamma$ represents a discrete group of orientation-preserving isometries. A hyperbolic $n$--orbifold is a \textit{manifold} when $\Gamma$ contain no elements of finite order. In \cite{Mar1}, Martin constructed a lower bound for $r_n$, the largest number such that every hyperbolic $n$--manifold contains a round ball of that radius; see also Friedland and Hersonsky \cite{FriHer}. From this one can compute, in each dimension, an explicit lower bound for the volume of a hyperbolic $n$--manifold.

The purpose of this paper is to give an explicit lower bound for the volume of a hyperbolic $n$--orbifold, again depending only on dimension. The result of this article is more general than what was achieved in the prequel \cite{Ade1}. Our work also significantly  improves upon the volume bounds of \cite{Ade1} and \cite{Mar1}, even though we consider a larger category of orbit spaces. 

We define a \textit{Riemannian submersion} $\pi\co SO_o(n,1)/\Gamma \to \mathbb{H}^n/ \Gamma $, where $SO_o(n,1)$, the connected component of the identity in the Lie group $O(n,1)$, is isomorphic to the full group of orientation-preserving isometries of $\mathbb{H}^n$. The study of the volume of a hyperbolic orbifold is thereby reduced to the study of the covolume of a lattice in a Lie group.

 In \cite{Wa1}, Wang showed that the covolume of a lattice in a semisimple Lie group that contains no compact factor can be bounded below by the volume of ball with a radius that depends only on the group itself. We estimate the sectional curvature of $SO_o(n,1)$ and apply a comparison theorem due to Gunther (see e.g. \cite{Gall}), to produce a lower bound for $\Vol[SO_o(n,1)/\Gamma]$. The following theorem gives our main result.

\begin{theorem}\label{Bound} The volume of a hyperbolic $n$--orbifold is bounded below by $\mathcal B(n)$, an explicit constant depending only on dimension, given by
\[\mathcal B(n)=\frac{2^{[\frac{6-n}{4}]}   \pi^{[\frac{n}{4}]}(n-2)!(n-4)!\cdots 1}{(2+9n)^{[\frac{n^2+n}{4}]}\Gamma(\frac{n^2+n}{4})}\int_{0}^{\min[0.08\sqrt{2+9n},\pi]}\sin^{\frac{n^2+n-2}{2}}\rho\,\,\,d\rho.\]
\end{theorem}

\begin{remark} The equation of Theorem~\ref{Bound} can be refined for $n=2,3$ to give a slightly better estimate. We describe these cases at the end of Section~\ref{last}.\end{remark}

The next section describes a canonical metric for $SO_o(n,1)$. Section 2 outlines Wang's crucial result. In the third section, we derive the curvature formulas for a canonical metric of a semisimple Lie group. These formulas are then used to construct an upper bound for the sectional curvatures of $SO_o(n,1)$. 

We prove Theorem~\ref{Bound} in Section~\ref{last}. From this formula, we get a lower bound of $2.46\times10^{-7}$ for hyperbolic $3$--orbifolds, $2.93\times10^{-13}$ for $4$--orbifolds and $2\times 10^{-20}$ for $5$--orbifolds. 

For comparison, the fifth section lists several results on hyperbolic volume. Sharp volume bounds for hyperbolic orbifolds are known for dimensions 2 and 3. The hyperbolic 2--orbifold of minimum volume was identified by Siegel \cite{CLS} in a theorem closely related to a result on birational transformations of an algebraic curve due to Hurwitz \cite{Hur}. The analogous result for dimension 3 was proved by Gehring and Martin \cite{GM1}. A hyperbolic orbifold is: a \textit{manifold} when $\Gamma$ does not contain elliptic elements; \textit{cusped} when $\Gamma$ does contain parabolic elements; \textit{arithmetic} when $\Gamma$ can be derived by a specific number-theoretic construction (see e.g. \cite{Belo}). What has been established for higher dimensions relate to these categories and their various intersections.

\bigskip

Intimately linked with hyperbolic volume is the size of symmetry groups of hyperbolic manifolds. Specifically, any bound in one category immediately produces a bound in the other. The quotient of a hyperbolic manifold $M$ by its group of orientation-preserving isometries is an orientable hyperbolic orbifold (as long as $\pi_1(M)$ is not virtually abelian, in which case $\Vol[M]$ is infinite). The following corollary is a direct analogue of Hurwitz's formula for groups acting on surfaces.

\begin{corollary}\label{cor1} Let $M$ be an orientable hyperbolic $n$--manifold. Let $H$ be a group of orientation-preserving isometries of $M$. Then 
\[
|H|\leq\frac{\Vol[M]}{\mathcal B(n)}.
\]
\end{corollary}

The Mostow--Prasad rigidity theorem \cite{mostow}, \cite{prasad} implies that the group of isometries of a finite volume hyperbolic $n$--manifold can be identified with $\Out(\pi_1(M))$. Hence, we have the following `topological' version of Corollary~\ref{cor1}.

\begin{corollary}\label{cor2} Let $M$ be a finite volume orientable hyperbolic $n$--manifold. Let $H$ be a subgroup of $\Out(\pi_1(M))$. Then 
\[
|H|\leq\frac{2\Vol[M]}{\mathcal B(n)}.
\]
\end{corollary}

\sect{The Canonical Metric of $SO_o(n,1)$}\label{canmet}

 Let $G$ be a Lie group and $\mathfrak g$ its Lie algebra. For $X \in \mathfrak{g}$, the \textit{adjoint action} of $X$ is the $\mathfrak g$-endomorphism defined by the Lie bracket \[\ad X(Y):=[X,Y].\]
 
 The \textit{Killing form} on $\mathfrak{g}$ is a symmetric bilinear form given by \[B(X,Y):=\tr(\ad X \circ \ad Y).\]  
 
 We note here that for all $X \in \mathfrak{g}$, $\ad X$ is skew symmetric with respect to $B$; i.e.,
  \be  \label{B-skew}
B([X, Y], Z) = -B(Y, [X,Z]).
\ee

A Lie group $G$ is called \textit{semisimple} if the Killing form associated to its Lie algebra is nondegenerate. In this case, there exists a \textit{Cartan decomposition}  $\mathfrak{g}=\mathfrak k \oplus \mathfrak{p}$ such that $B\, | \, \mathfrak k$ is negative definite and $B\, | \, \mathfrak p$ is positive definite,  with bracket laws 
\be  [\mathfrak k, \mathfrak k] \subset \mathfrak k, \  \ [\mathfrak k, \mathfrak p] \subset \mathfrak p, \  \ [\mathfrak p, \mathfrak p] \subset \mathfrak k.   \label{brackets}
\ee

A positive definite inner product on $\mathfrak{g}$ is defined by putting 

\[ \langle X,Y\rangle := \left\{ \begin{array}{lll} B(X,Y) & \mbox{for } X,Y\in\mathfrak p \\
-B(X,Y) & \mbox{for } X,Y\in\mathfrak k \\ 0 & \mbox{otherwise}. \end{array} \right.\]

Let $e$ denote the identity element of $G$. We identify $\mathfrak{g}$ with $T_eG$, the tangent space of $G$ at the identity, and extend $\langle X,Y\rangle$ to a left invariant Riemannian metric over $G$ by left translation. This metric will be referred to as a \textit{canonical metric} for $G$. When the choice of Cartan decomposition is clear, we denote the associated canonical metric by $g$ and the induced distance function on $G$ by $\rho$.

Let $K$ denote the maximal compact subgroup of $G$ with Lie algebra $\mathfrak k$. Important in what follows is that the restriction of $\langle X,Y\rangle$ to $\mathfrak p\simeq T_{eK} G/K$ induces a Riemannian metric on the quotient space, as well. In Definition~\ref{mainmetric}, the canonical metric on a specific Lie group $G$ is scaled in order to secure desired curvature properties for $G/K$.

\bigskip

Denote by $GL(n,\mathbb{R})$ the group of real nonsingular $n$-by-$n$ matrices. The \textit{Lorentz group} $O(n,1)$ is defined by \[O(n,1):=\{A\in GL(n+1,\mathbb{R}):JA^TJ=A^{-1}\}, \mbox{ where } J=\left(\begin{smallmatrix}
1 &   &   &   &   \\  & 1 &   &   &   \\  &   & \ddots &   &   \\  &   &   & 1 &   \\  &   &   &   & -1
\end{smallmatrix}\right).\] 

The Lorentz group is a \textit{matrix Lie group}; it is a differentiable manifold where matrix multiplication is compatible with the smooth structure. The \textit{positive special Lorentz group} $SO_o(n,1)$ is the identity component of $O(n,1)$. It consists of the elements of $O(n,1)$ that have determinant $1$ and a positive $(n+1,n+1)$ coordinate.

The Lie algebra of any matrix Lie group $G$ is the set of matrices $X$ such that  $e^{tX}\in G$, for all real numbers $t$. Denote by $\mathfrak{so}(n,1)$ the Lie algebra of $SO_o(n,1)$. Then 
\begin{eqnarray*}
X \in \mathfrak{so}(n,1)&\Rightarrow&e^{tX}\in SO_o(n,1)\\
&\Rightarrow& J(e^{tX})^TJ=(e^{tX})^{-1}\\
&\Rightarrow& Je^{tX^T}J=e^{-tX}\\
&\Rightarrow& e^{tJX^TJ}=e^{-tX}\\
&\Rightarrow& JX^TJ=-X.\\
\end{eqnarray*}

Let $X=(a_{ij})$ be an $n+1$-by-$n+1$ matrix. If $JX^TJ=-X$, then $X$ has the form
\[\left(\begin{array}{ccccc}0 & a_{12} & a_{13} & \hdots & a_{1,n+1} \\-a_{12} & 0 & a_{23} & \hdots & a_{2,n+1} \\-a_{13} & -a_{23} & 0 &   &  \\\vdots &  &   & \ddots &   \\a_{1,n+1} & a_{2,n+1} &   &   & 0\end{array}\right).\]

For each $n$, let $e_{ij}$ represent the $n+1$-by-$n+1$ matrix with $1$ in the $ij$-position and $0$ everywhere else. Let $\alpha_{ij}=( e_{ij}-e_{ji})$ and $\sigma_{ij}=( e_{ij}+e_{ji}).$

\begin{definition} \label{basiselt} The \textbf{standard basis} for  $\mathfrak{so}(n,1)$, denoted by $\mathfrak{B}$, consists of the following set of  $n(n+1)/2$ matrices:
\begin{align*}
\alpha_{12} &&\alpha_{13} &&\alpha_{14} &&\hdots &&\alpha_{1n} &&\sigma_{1,n+1}\\
&&\alpha_{23} &&\alpha_{24} &&\hdots &&\alpha_{2n} &&\sigma_{2,n+1}\\
&&&&\alpha_{34} &&\hdots &&\alpha_{3n} &&\sigma_{3,n+1}\\
&&&&&&\ddots&&&&\\
&&&&&&&&\alpha_{n-1,n} &&\sigma_{n-1,n+1}\\
&&&&&&&&&&\sigma_{n,n+1}\\
\end{align*}
\end{definition}

The Lie bracket of a matrix Lie algebra is determined by matrix operations.  \[ [X,Y] := XY-YX.\]
 The following proposition describes the Lie bracket of  $\mathfrak{so}(n,1)$.
\begin{prop}  \label{multtable}
For $1\le i < j\le n, 1\le k < l\le n,$
\ba
[\alpha_{ij}, \alpha_{kl}] & = & \delta_{jk} \alpha_{il} + \delta_{jl}\alpha_{ki} + \delta_{il} \alpha_{jk}+ \delta_{ki}\alpha_{lj}  \\
& = & \left\{ \begin{array}{ll}
\alpha_{il} & \mbox{if} \ j=k \\
\alpha_{ki} & \mbox{if} \ j=l \\ 
\alpha_{jk} & \mbox{if} \ i=l \\
\alpha_{lj} & \mbox{if} \ i=k \\
0 & \mbox{otherwise}
\end{array} \right.    \label{bracket-al-al}
\ea

\ba
[\alpha_{ij}, \sigma_{k, n+1}] & =&  \delta_{kj}  \sigma_{i,n+1} - \delta_{ik} \sigma_{j, n+1}   \label{bracket-alpha-sigma}  \\
& = & \left\{ \begin{array}{ll}
\sigma_{i,n+1} & \mbox{if} \ k=j \\
-\sigma_{j,n+1} & \mbox{if} \ i=k \\ 
0 & \mbox{otherwise}
\end{array} \right.  
\ea

\be\label{sigref}
[ \sigma_{i,n+1}, \sigma_{j,n+1}]  =  \alpha_{ij}.
\ee
\end{prop}
\begin{proof}The proof of the first equation is given here. The proofs of the remaining identities are similar. 

By the definition of $\alpha_{ij}$ and the fact that  $e_{ij}e_{kl} = \delta_{jk}e_{il}$,
\begin{eqnarray*}
[\alpha_{ij},\alpha_{kl}] &=& [e_{ij}-e_{ji},  e_{kl}-e_{lk}]  \\
& = & (e_{ij}-e_{ji})(e_{kl}-e_{lk})-(e_{kl}-e_{lk})(e_{ij}-e_{ji})\\
&=& e_{ij}e_{kl}-e_{ij}e_{lk}-e_{ji}e_{kl}+e_{ji}e_{lk}-e_{kl}e_{ij}+e_{kl}e_{ji}+e_{lk}e_{ij}-e_{lk}e_{ji}\\
&=& \delta_{jk} (e_{il} - e_{li}) + \delta_{jl} (e_{ki} -e_{ik} ) + \delta_{il} (e_{jk}- e_{kj}) + \delta_{ki} (e_{lj} - e_{jl})  \\
&=& \delta_{jk} \alpha_{il} + \delta_{jl}\alpha_{ki} + \delta_{il} \alpha_{jk}+ \delta_{ki}\alpha_{lj} .
\proved\end{eqnarray*}\end{proof}

Proposition~\ref{multtable} illustrates a Cartan decomposition $\mathfrak{so}(n,1) = \mathfrak{k} \oplus \mathfrak{p}$, where \be\mathfrak{k} =\spn\{ \alpha_{ij}, 1\le i <j \le n\} \mbox{ and } \mathfrak{p}= \spn\{  \sigma_{i,n+1}, 1\le i\le n\}.\label{CD}\ee 
We note here that $\mathfrak{k}=\mathfrak{so}(n)$, the Lie algebra of the Lie group $SO(n)$. In turn, $SO(n)$ is a maximal compact subgroup of $SO_o(n,1)$.

In this article, \textit{the canonical metric for $SO_o(n,1)$} refers to the canonical metric induced by the Cartan decomposition of (\ref{CD}). It is denoted by $g$. The following lemma and corollary give a description of the metric $g$.

\begin{lemma} Let $X,Y\in\mathfrak B$. Then

\begin{displaymath}
\lp X,Y \rp =\left\{\begin{array}{ccr}
2n-2&\rm{if     }  &X=Y \\
0& &\rm{otherwise.}
\end{array}
\right.
\end{displaymath}
\end{lemma}

\begin{proof} 

The proof follows from a close study of Proposition~\ref{multtable}. 

The set $\mathfrak B$ is closed under the Lie bracket (modulo sign). Therefore, for any $X\in \mathfrak B$ the entries of $\ad X$ are all $0,1$ or $-1$ and each column has at most one non-zero entry. Since bracket multiplication is determined by index, each row also has at most one non-zero entry. Furthermore, two standard basis elements have a non-zero Lie bracket if and only if they share exactly one index number. So if $X$ has index $ij$, $\ad X$ has exactly \[(n+1-i)+(j-1)+(n+1-j)+(i-1)-1-1=2n-2\] non zero entries.

Now assume $X=\alpha_{ij}$. For all $Y\in\mathfrak B, [X,Y]=Z\Rightarrow [X,Z]=-Y.$ This implies that the $hg$ entry of $\ad X$ is the negative of the $gh$ entry.

By definition,
\begin{eqnarray*}
\lp \alpha_{ij},\alpha_{ij} \rp&=& -B(\alpha_{ij},\alpha_{ij})\\
&=& -\tr(\ad \alpha_{ij}\circ\ad \alpha_{ij}).
\end{eqnarray*}
The $h$th diagonal entry of $\ad \alpha_{ij}\circ\ad \alpha_{ij}$ is the dot product of the $h$th row of $\ad \alpha_{ij}$ with the $h$th column of $\ad \alpha_{ij}$. If the only non-zero entry in the $h$th row of $\ad \alpha_{ij}$ is a $1$ (resp. $-1$) in the $hg$-position then the only non-zero entry  in the $h$th column of $\ad \alpha_{ij}$ is a $-1$ (resp. $1$) in the $gh$-position. Hence, the $h$th diagonal entry of $\ad \alpha_{ij}\circ\ad \alpha_{ij}$ is $-1$. Thus, \[\lp\alpha_{ij},\alpha_{ij}\rp=-\left(\underbrace{-1+-1+\cdots+-1}_{2n-2 \text{    times}}\right)=2n-2.\]Similarly, $\lp\sigma_{ij},\sigma_{ij}\rp=2n-2$. 

Let $X,Y\in\mathfrak B$, with $X\neq \pm Y$. If $\ad X$ has a nonzero entry in the $hg$-position then the bracket of $X$ with the $h$th basis element is sent to the $g$th basis element. That is, there exists $V, W\in\mathfrak B$ such that
\[[X,V]=\pm W.\]
If, in addition, $\ad Y$ has a nonzero entry in the $gh$-position, we may write
\[[Y,W]=\pm V.\]
Again, note that the Lie bracket of basis elements is determined by index. This forces \[X=\pm Y\] and we have a contradiction. Thus, all the diagonal entries of $\ad X \circ\ad Y$ are equal to zero. Therefore $\lp X,Y\rp=0$.
\end{proof}

\begin{corollary}  \label{metric}
The matrix representation for $g$, the canonical metric for $SO_o(n,1)$, is the square $n(n+1)/2$ diagonal matrix \[\left(\begin{array}{cccc}2n-2 &   &   &   \\  & 2n-2 &   &   \\  &   & \ddots &   \\  &   &   & 2n-2\end{array}\right).\]
\end{corollary}

We will be interested in the metric that induces constant sectional curvature $-1$ on the quotient space $SO_o(n,1)/SO(n)$. To this end, we scale the metric $g$ by the factor $\frac{1}{2n-2}$.

\begin{definition}\label{mainmetric} Let $g$ be the canonical metric for $SO_o(n,1)$. The metric $\tilde g$ on $SO_o(n,1)$ is defined by \[\tilde g:=\frac{1}{2n-2}g.\]\end{definition}

\sect{Discrete Subgroups of Semisimple Lie Groups}

For any Lie group, a theorem of Zassenhaus \cite{Zass} guarantees the existence of a neighborhood $U$ of the identity such that the subgroup generated by any subset of $U$ is either non-discrete or nilpotent. Such a neighborhood is called a \textit{Zassenhaus} neighborhood. 

Kazhdan and Margulis \cite{KM} proved that if $G$ is a semisimple Lie group without compact factor it contains a Zassenhaus neighborhood $U$ such that, for any discrete subgroup $\Gamma$ of $G$, there exists $g\in G$ with the property that $g\Gamma g^{-1}\cap U=\{e\}$. This implies that the fundamental domain for any lattice in $G$ has a definite size.

In  \cite{Wa1}, H. C. Wang undertook a quantitative study of a Zassenhaus neighborhood for a semisimple Lie group $G$, with respect to a canonical metric. Wang found a value $R_G$ such that a metric ball in $G$ centered at the identity with radius $R_G$ satisfied the conclusion of the Kazhdan-Margulis theorem.

Recall the definitions and notations of Section~\ref{canmet}. Again, let $G$ be a Lie group, $\mathfrak g$ its Lie algebra, $\mathfrak g=\mathfrak k \oplus\mathfrak p$ a Cartan decomposition and $\langle\cdot,\cdot\rangle$ the associated inner product. Define a norm on $\mathfrak g$ by $\norm X:=\langle X, X \rangle^{1/2}$. For each $\mathfrak g$-endomorphism $f$, let \[N(f):=\sup\{||f(X)||:X\in\mathfrak g, ||X||=1\}.\] Furthermore, let 
\begin{gather*}
C_1:=\sup\{N(\ad X):X\in\mathfrak p, ||X||=1\}\\
\tag*{\text{and}} C_2:=\sup\{N(\ad X):X\in\mathfrak k, ||X||=1\}.
\end{gather*}

The number $R_G$ is defined to be the least positive zero of the real-valued function \be\label{Fval}F(t)=\exp C_1t-1+2\sin C_2t-\frac{C_1t}{\exp C_1t-1}~.\ee

The following theorem (Theorem 3.2 in \cite{Wa1}) demonstrates the role of the value $R_G$ in the construction of a Zassenhaus neighborhood for a semisimple Lie group.

\begin{theorem}[Wang]\label{wang1} Let $G$ be a semisimple Lie group. Let $e\in G$ denote the identity. Then for any discrete subgroup $\Gamma$ of G, the set \[\Theta=\{g\in\Gamma:\rho(e,g)\leq R_G\}\] generates a nilpotent group.
\end{theorem}

Now, let $\mathfrak{g}_{\pi}$ be the totality of elements $X$ in $\mathfrak{g}$ such that the imaginary parts of all the eigenvalues of $\ad X$ lie in the open interval $(-\pi,\pi)$ and let $G_{\pi}=\{\expo X:X\in \mathfrak{g}_{\pi}\}$. In an earlier work \cite{Wa3}, Wang had proved that the restriction of the exponential map to $\mathfrak{g}_{\pi}$ is injective. Hence, the following proposition (Proposition 5.1 in \cite{Wa1}) establishes the fact that $R_G$ is less than the injectivity radius of $G$. 

\begin{prop}[Wang]\label{injrad} Let $G$ be a semisimple Lie group. Then the closed ball \[B_G=\{x\in G:\rho(e,x)\leq R_G\}\]is contained in $G_{\pi}$.
\end{prop}

We now give Wang's quantitative version of the theorem of Kazhdan-Margulis (Theorem 5.2 in \cite{Wa1}). It shows that the volume of the fundamental domain of $\Gamma$ is larger than the volume of a $\rho$-ball with radius $R_G/2$.

\begin{theorem}[Wang]\label{wang2} Let $G$ be a semisimple Lie group without compact factor and $B_G=\{x\in G:\rho(e,x)\leq R_G\}$. Then for any discrete subgroup $\Gamma$ of G, there exists $g\in G$ such that $B_G\cap g \Gamma g^{-1} =\{e\}$.\end{theorem}

The appendix to \cite{Wa1} includes a table of the constants $C_1$ and $C_2$ for noncompact and nonexceptional Lie groups. For $SO_o(n,1)$, $n\geq 4$, with respect to the scaled canonical metric $\tilde g$ (Definition~\ref{mainmetric}), we have  \be\label{Cval} C_1=1 \mbox{ and } C_2=\sqrt 2.\ee
Therefore, by (\ref{Fval}) and (\ref{Cval}), \be\label{rg} R_G=228/1000 \text{\,\,\,when\,\,\,}G=SO_o(n,1), n\geq 4.\ee
When $n=2,3$, \be\label{C3val} C_1=C_2=1.\ee This gives
\be\label{rg3} R_G=277/1000 \text{\,\,\,when\,\,\,}G=SO_o(n,1),  n=2,3.\ee

\sect{The Sectional Curvatures of $SO_o(n,1)$}

In this section, we construct an upper bound on the sectional curvatures of $SO_o(n,1)$. As a first step, we derive the curvature formulas for a canonical metric of a semisimple Lie group. These formulas are of independent interest as we could not find them in the literature.

A \textit{connection} $\nabla$ on the tangent bundle of a manifold can be expressed in terms of a left invariant metric by the \textit{Koszul formula}. For any left invariant vector fields $X,Y,Z,W$, we have
\be
\lp \nabla _XY , Z \rp =  \frac 12\left\{ \lp [X,Y], Z\rp   - \lp Y, [X, Z] \rp -  \lp X, [Y,Z]\rp  \right\}.  \label{conn-l}
\ee
The \textit{curvature tensor} of a connection $\nabla$ is defined by
\be R(U,V)X=\nabla_U\nabla_VX-\nabla_V\nabla_UX-\nabla_{[U,V]}X.  \label{curvature-tensor} \ee

When a Lie group $G$ is semisimple and compact, the canonical metric is the negative of the Killing form and induces a biinvariant metric on $G$. The connection and curvature can be described in terms of the Lie bracket in a simple way (see e.g. \cite[Cor. 3.19]{Cheeger-Ebin}).
\ba
\nabla _XY & = &  \frac 12 [X,Y], \\
\lp R(X,Y)Y,X \rp & = & \frac 14 \|[X,Y]\|^2.
\ea

When $G$ is semisimple and noncompact, a canonical metric is biinvariant only when restricted to $K$, the maximal compact subgroup of $G$ with Lie algebra $\mathfrak k$. The connection and curvature formulas for this case are given below.

We will treat vector fields from $ \mathfrak k$  and $ \mathfrak p$ separately.   From here on, $U,V, W, W' \in \mathfrak k$ and $X,Y, Z, Z' \in \mathfrak p$ denote left invariant vector fields. 

\begin{lemma}\label{totgeo}
With respect to the canonical metric the subgroup $K$ is totally geodesic in $G$. 
\end{lemma}
\begin{proof}
Since the canonical metric restricted to $K$ is biinvariant \[\lp \nabla _U U, V\rp = 0.\] By (\ref{brackets}) and (\ref{conn-l}),   \[\lp \nabla _U U, X\rp =  - \lp U, [U,X] \rp =  0.\proved\]\end{proof}

Now we compute the connections.

\begin{lemma}  \label{connections}
\ba
\nabla _UV & = &  \frac 12 [U,V],  \\
 \nabla _XY & = &  \frac 12 [X,Y], \\
  \nabla _UX   =   \frac 32 [U,X], & &   \nabla _X U  =  - \frac 12 [X,U].
  \ea
  \end{lemma}
 \begin{proof}
 The first equation follows from Lemma~\ref{totgeo}. We will derive the last two equations, the proof of the second equation is similar.
 
 Again by (\ref{brackets}) and (\ref{conn-l}), \[
  \lp \nabla _UX, V \rp = \frac 12 \lp [U,X], V \rp - \frac 12 \lp X, [U,V] \rp - \frac 12 \lp U, [X,V] \rp=0
  \]
  and 
  \[
  \lp \frac 32 [U,X],V\rp =0.
  \]
  Thus,
  \be
    \lp \nabla _UX, V \rp = \lp \frac 32 [U,X],V\rp, \mbox{ for all }~ V\in\mathfrak k.\label{UXV}
    \ee
  
  Similarly, by (\ref{B-skew}) and (\ref{conn-l}),
  \[
  \lp \nabla _UX, Y \rp = \frac 12 \lp [U,X], Y \rp - \frac 12 \lp X, [U,Y] \rp - \frac 12 \lp U, [X,Y] \rp
  \]
  and
 \[\lp X, [U,Y] \rp= - \lp [U,X], Y \rp\] and
  \[ \lp U, [X,Y] \rp = -  B( U, [X,Y]) = B( [X,U], Y ) = \lp [X,U], Y \rp.\]
 Hence,
 \be
  \lp \nabla _UX, Y \rp = \frac 12 \lp [U,X], Y \rp + \frac 12 \lp [U,X], Y \rp + \frac 12 \lp [U,X], Y \rp=\lp \frac 32 [U,X], Y \rp, \mbox{for all }Y\in\mathfrak p. \label{UXY}
  \ee 
 
 From  (\ref{UXV}) and (\ref{UXY}), we have
 \[\nabla _UX   =   \frac 32 [U,X].\]
 
 Finally, 
 \[\nabla _X U =  \nabla _UX + [X,U] =  - \frac 12 [X,U].\proved\] \end{proof}
  
The following proposition gives the corresponding curvature formulas.
\begin{prop}  \label{curvatures}
\ba
R(U,V)W & = & \frac 14 [ [V,U], W],   \label{R-UVW}  \\
 R(X,Y)Z & = & - \frac 74 [ [X,Y],  Z ],    \label{R-XYZ}\\
 R(U,X)Y & = & \frac 14 [[X,U], Y]  - \frac 12 [[Y,U], X] , \label{R-UXY}  \\
 R(X,Y)V & = & \frac 34 [X, [V,Y]] + \frac 34 [Y, [X,V]].  \label{R-XYV}
\ea
In particular,
\ba
\lp R(U,V)W, X \rp & = & 0,   \label{R-mix1}\\
\lp R(X,Y)Z, U\rp & = &0,  \label{R-mix2}  \\
\lp R(U,V)V,U\rp & =  & \frac 14 \|[U,  V]\|^2,  \label{R-UV}\\
\lp R(X,Y)Y,X \rp & = & -  \frac 74  \|[X,  Y]\|^2, \label{R_XY}\\
\lp R(U,X)X,U\rp & =  & \frac 14   \|[U,  X]\|^2.   \label{R-UX}
\ea
\end{prop}
\begin{proof}
We prove (\ref{R-XYZ}). The proofs of the remaining equations are similar.

By (\ref{brackets}), (\ref{curvature-tensor}) and Lemma~\ref{connections}, 
\ban
 R(X,Y)Z &  = & \frac 12 \left(   \nabla _X [Y,Z]  -  \nabla _Y [X,Z] - 3 [[X,Y], Z]  \right)  \\ 
 & = &  \frac 12 \left( -\frac 12[X, [Y, Z] ] +  \frac 12 \lp [Y, [X,Z] ]  - 3 [[X,Y], Z]  \right).
\ean
Therefore, by the Jacobi identity,
\[ [X,[Y,Z]]+[Y,[Z,X]]+[Z,[X,Y]]=0 \mbox { for all } X,Y,Z \in \mathfrak g,\]
\ban
 R(X,Y)Z   & = &   \frac 12 \left(  \frac 12   [Z,  [X, Y]]  - 3 [[X,Y], Z]  \right)  \\
 & = & -  \frac 74 [[X,Y], Z] .\proved\ean\end{proof}
 
 \bigskip
  
 For a Lie group $G$, with Lie algebra $\mathfrak g$ and $X,Y \in \mathfrak g$, the \textit{sectional curvature} of the planes spanned by $X$ and $Y$ is denoted and defined by \[K(X,Y)=\frac{\lp R(X,Y)Y,X \rp}{\|X\|^2\|Y\|^2-\lp X,Y\rp^2}.\] 
 
In the next two propositions, we develop our bound for the sectional curvatures of $SO_o(n,1)$. Recall the notation established in Section~\ref{canmet}.

\begin{prop}\label{seccurbas}
The sectional curvature of $SO_o(n,1)$ with respect to the metric $\tilde g$ at the planes spanned by standard basis elements is bounded above by $\frac{1}{4}$. 
\end{prop}
\begin{proof}
Since $\alpha_{ij}, \alpha_{kl}$ are orthogonal,\[K(\alpha_{ij}, \alpha_{kl}) =\frac{ \langle R(\alpha_{ij},\alpha_{kl})\alpha_{kl},\alpha_{ij} \rangle }{\|\alpha_{ij}\|^2 \|\alpha_{kl}\|^2}.\]
By (\ref{R-UV}), Proposition~\ref{multtable} and Corollary~\ref{metric},
\be
K(\alpha_{ij}, \alpha_{kl})   =   \frac{ \|[\alpha_{ij},  \alpha_{kl}]\|^2}{4 \|\alpha_{ij}\|^2 \|\alpha_{kl}\|^2}  \le   \frac{ 1}{4}.
\ee
Similarly,
\be
K(\alpha_{ij},\sigma_{k,n+1})  =    \frac{ \|[\alpha_{ij},  \sigma_{k,n+1}]\|^2 }{ 4\|\alpha_{ij}\|^2 \|\sigma_{k,n+1}\|^2}   \le  \frac{ 1}{4}\ee
and 
\be\label{sigref2} K(\sigma_{i,n+1},\sigma_{j,n+1})   =  -  \frac {7 \|[\sigma_{i,n+1},  \sigma_{j,n+1}]\|^2}{4\|\sigma_{i,n+1}\|^2 \|\sigma_{j,n+1}\|^2} = \frac{-7}{4}.
\ee
\end{proof}

\begin{prop}\label{seccur} The sectional curvatures of $SO_o(n,1)$ with respect to $\tilde{g_0}$ are bounded above by 
\[
\frac 12 + 2 \frac{1}{4} +  2 \frac{6n}{4}  + 2 \frac{3n}{4}  =  \frac{2+ 9n}{2}. \]
\end{prop}

\begin{remark}\label{twothree} Using (\ref{C3val}) instead of (\ref{Cval}) in the proof of Proposition~\ref{seccur} gives a bound of $(3+18n)/4$ for $n=2,3$. In dimension 2, additional calculation reduces the bound to $1/4$.
\end{remark}

\begin{proof} Again with $U,V \in \mathfrak{k}\mbox{ and } X, Y \in \mathfrak{p}$, we have by (\ref{R-mix1}) and (\ref{R-mix2}), 
\begin{eqnarray*}
\langle R(X+U,Y+V)Y+V,X+U \rangle  & = & \langle R(X,Y)Y,X\rangle + \langle R(U,V)V,U\rangle + \langle R(U,Y)Y,U\rangle \\
&   +& \langle R(X,V)V,X \rangle  + 2\langle R(X,Y)V, U\rangle + 2\langle R (X,V)Y,U \rangle. 
\end{eqnarray*}
Assume that $\|U+X\|=1,~ \|V+Y\| = 1 \mbox{ and }\langle U+X, V+Y \rangle = 0$. 
Write \[U = \sum_{i<j} a_{ij} \alpha_{ij},~  V = \sum_{i<j} a_{ij}' \alpha_{ij},~ X = \sum_{i=1}^n b_i \sigma_{i, n+1}, ~ Y = \sum_{i=1}^n b_i'  \sigma_{i, n+1}.\]  
Note that 
\begin{equation}
\sum_{i<j} |a_{ij}|^2,\   \sum_{i<j} |a'_{ij}|^2, \ \sum_{i=1}^n  |b_i|^2,   \ \sum_{i=1}^n  |b'_i|^2 \le 1.
\end{equation}
 By (\ref{R-UVW}),
 \[
 R(U,V)V = \frac 14 [[V,U], V]  = -\frac 14  \ad V \circ \ad V (U).
 \]
 Hence, 
 \[
 \langle R(U,V)V,U\rangle \le \frac 14 C_2^2  = \frac 12.\]

Similarly, by (\ref{R-UXY}),
\[
R(U,Y)Y   =  - \frac 14 [[Y,U],Y] = \frac 14 \ad Y\circ \ad Y (U)
\]
and
\[
\langle R(U,Y)Y,U\rangle \le \frac 14 C_1^2 = \frac 14.
\]

By (\ref{R-XYV}),
\[
\langle R(X,Y)V, U\rangle = - \frac 34 \left( \lp [U,X], [V,Y] \rp + \lp [V,X], [U,Y] \rp \right).
\]
Now 
\begin{eqnarray*} 
\| [U,Y] \|^2   & = &    \| \left[\sum_{i<j} a_{ij} \alpha_{ij},  \sum_{k} b'_{k} \sigma_{k,n+1}\right]\|^2  \\
& = &    \|  \sum_k \left( \sum_i a_{ki} b_i' \right) \sigma_{k, n+1} \|^2 \\
& = &      \sum_k \left( \sum_i a_{ki} b_i' \right)^2  \le  n.
 \end{eqnarray*}
Hence, \[
\langle R(X,Y)V, U\rangle \le  \frac{6}{4} \cdot n.
\]
Similarly,  by (\ref{R-UXY}),
\[
\langle R (X,V)Y,U \rangle \le \frac{3}{4} \cdot  n. 
\proved\]
\end{proof}

\section{Volumes of Hyperbolic $n$--Orbifolds}
\label{last}

Let $(M,g)$ and $(N,h)$ be Riemannian manifolds and $q:M\rightarrow N$ a surjective submersion. For each point $x\in M$ the tangent space $T_xM$ decomposes into the orthogonal direct sum \[T_xM=(\Ker dq)^{\bot}_{x}\oplus (\Ker dq)_{x}.\]The map $q$ is said to be a \textit{Riemannian submersion} if \[g(X,Y)=h(dq X,dq Y)\text{\,\,\,whenever\,\,\,}X,Y\in(\Ker dq)^{\bot}_{x}\text{\,\,\,for some\,\,\,} x\in M.\]

\begin{lemma}\label{resub}
Let $K \ra M \stackrel{q}{\ra} N$ denote a fiber bundle where $q$ is a Riemannian submersion and  $K$ is a compact and totally geodesic submanifold  of $M$. Then for any subset $Z \subset N$, \[\vol [q^{-1}(Z)] =  \vol[Z] \cdot \vol[K].\]
\end{lemma}
\begin{proof}
Since $K$ is totally geodesic, the fibers of $q$ are isometric to each other.  Therefore,
\[
d \vol_M = d \vol_N \cdot d \vol_K.
\]
Hence, $\vol[q^{-1}(Z)] = \int_{q^{-1}(Z)} d \vol_M  = \int_{q^{-1}(Z)}  d \vol_K  \, d \vol_N = \vol (K) \cdot \vol(Z)$.
 \end{proof}
 
 Let \[\pi\co SO_o(n,1) \to SO_o(n,1)/SO(n)\] be the quotient map. Recalling the definitions and notations of Section~\ref{canmet}, equip $SO_o(n,1)$ with the scaled canonical metric $\tilde g$. Furthermore, assign to the quotient the metric induced by the restriction of $\tilde g$ to $\mathfrak p\subset\mathfrak{so}(n,1)$. The map $\pi$ is then a Riemannian submersion.
 
 O'Neill's formula (see e.g. \cite[Page 127]{Gall}) relates the sectional curvature of the base space of a Riemannian submersion, $K_b$, with that of the total space, $K_t$. Let $X,Y\in\mathfrak p$ represent orthonormal vector fields on $SO_o(n,1)/SO(n)$ as well as their horizontal lifts. O'Neill's formula, applied to $\pi$, gives \[K_b(X, Y)=K_t(X, Y)+\frac 3 4|[X, Y]^v|^2.\] Here, $Z^v$ denotes the vertical component of $Z$. 
  
 From (\ref{CD}) and (\ref{R_XY}), we then get \be K_b(X,Y)=-  \frac 74  \|[X,  Y]\|^2+\frac 3 4 \|[X,  Y]\|^2=- \|[X,  Y]\|^2,\label{scbrs}\ee
 where \[X = \sum_{i=1}^n a_i \sigma_{i, n+1},~Y = \sum_{i=1}^n b_i  \sigma_{i, n+1}. \]
 By Proposition~\ref{multtable}, Corollary~\ref{metric} and Definition~\ref{mainmetric}, \[\|[X,  Y]\|^2=\sum_{i<j} (a_ib_j-a_jb_i)^2.\]
 Since $\sum a_i^2=1,\sum b_i^2=1$ and $\sum a_ib_i=0$, we have 
 \begin{eqnarray*}2\sum_{i<j} (a_ib_j-a_jb_i)^2&=&\sum_{ij} (a_ib_j-a_jb_i)^2\\
 &=&\sum_{ij} a_i^2b_j^2 + \sum_{ij} a_j^2b_i^2 - 2\sum_{ij} a_ib_ja_jb_i \\ &=& 2. \end{eqnarray*}
 
 It follows that the quotient space $SO_o(n,1)/SO(n)$, with respect to the restriction of the scaled canonical metric, has constant sectional curvature \[K_b(X,Y)=-1.\]
 It can therefore be identified with hyperbolic space, $\mathbb{H}^n$.
 
 For a discrete group $\Gamma < SO_o(n,1)$, let $Q$ be the hyperbolic $n$--orbifold defined by the quotient $\mathbb{H}^n / \Gamma$. The map $\pi$ induces another Riemannian submersion 
 \[\pi^{\prime}\co SO_o(n,1)/\Gamma  \to  Q.\]
 The fibers of $\pi^{\prime}$ on the smooth points of $Q$ are totally geodesic embedded copies of $SO(n)$. By Lemma~\ref{resub}, we have 
  \be\Vol[SO_o(n,1)/\Gamma]=\Vol[Q]\cdot\Vol[SO(n)]\label{frac}.\ee
  
  \bigskip
  
  Denote by $V(d,k,r)$ the volume of a ball of radius $r$ in the  complete simply connected Riemannian manifold of dimension $d$ with constant curvature $k$. A proof of the following comparison theorem can be found in \cite[Theorem 3.101]{Gall}.

\begin{theorem}[Gunther]\label{gun}Let $M$ be a complete Riemannian manifold of dimension $d$. For $m\in M$, let $B_{m}(r)$ be a ball which does not meet the cut-locus of $m$.

If the sectional curvatures of $M$ are bounded above by a constant $b$, then
\[\Vol[B_{m}(r)]\geq V(d,b,r).\]
\end{theorem}

\begin{prop}\label{numer} Let $\Gamma$ be a discrete subgroup of $SO_o(n,1)$. Then \[\Vol[SO_o(n,1)/\Gamma] \geq V(d_0,k_0,r_0),\] where $d_0=\dfrac{n^2+n}{2}$, $k_0=\dfrac{2+9n}{2}$ and $r_0=0.114$.
\end{prop}
\begin{proof} By Proposition~\ref{injrad}, Theorem~\ref{wang2} and (\ref{rg}), the volume of a fundamental domain of $\Gamma$ in $SO_o(n,1)$ is greater than the volume of a ball of radius $0.114$. 

From Definition~\ref{basiselt}, the dimension of $SO_o(n,1)$ is $(n^2+n)/2$. By Proposition~\ref{seccur} the sectional curvatures of $SO_o(n,1)$ are bounded above by $(2+9n)/2$. The Proposition then follows from Theorem~\ref{gun}.
\end{proof}

\bigskip 
 
 In \cite[Page 399]{Gilmore}, the volumes of the classical compact groups are given explicitly. For the special orthogonal group, the volume with respect to the metric $\tilde g$ is given by
 
 \be\Vol[SO(n)]=\frac{2^{[\frac{n^2+2n-2}{4}]}\pi^{[\frac{n^2}{4}]}}{(n-2)!(n-4)!\cdots1}\label{SON}\,\,.\ee
 
 \bigskip
We now prove Theorem~\ref{Bound}, which for convenience is restated below.

\medskip
{\bf Theorem \ref{Bound}}\qua{\sl The volume of a hyperbolic $n$--orbifold is bounded below by $\mathcal B(n)$, an explicit constant depending only on dimension, given by
\[\mathcal B(n)=\frac{2^{[\frac{6-n}{4}]}   \pi^{[\frac{n}{4}]}(n-2)!(n-4)!\cdots 1}{(2+9n)^{[\frac{n^2+n}{4}]}\Gamma(\frac{n^2+n}{4})}\int_{0}^{\min[0.08\sqrt{2+9n},\pi]}\sin^{\frac{n^2+n-2}{2}}\rho\,\,\,d\rho.\]}

\begin{proof} For $k>0$, the complete simply connected Riemannian manifold with constant curvature $k$ is the sphere of radius $k^{-1/2}$. By explicit computation we have \[V(d,k,r)= \frac{2(\pi/k)^{d/2}}{\Gamma(d/2)}\int_{0}^{\min\left[rk^{1/2},\pi\right]}\sin^{d-1}\rho\,\,\,d\rho.\]

The proof now follows from Proposition~\ref{numer}, (\ref{frac}) and (\ref{SON}).\end{proof}

In light of (\ref{rg3}) and Remark~\ref{twothree}, we can restate Proposition~\ref{numer} for $n=2,3$. In both cases, we have $r_0=0.1385$. The value for $k_0$ is taken to be $0.25$ when $n=2$ and $14.25$ when $n=3$. By the proof of Theorem~\ref{Bound}, our improved bounds for $n=2,3$ are \begin{align*}\mathcal B(2) &= 1\times10^{-3},  \\ \mathcal B(3) &=   2.46\times10^{-7}.\end{align*}

\section{Volume Bounds}
\label{comparisons}

That the smallest hyperbolic 2--orbifold and 2--manifold have area, respectively, $\pi/21$ and $4\pi$ are classical results. In \cite{GM1}, it was proved that the smallest hyperbolic 3--orbifold has volume 0.03905.... It was shown in \cite{GabMeyMilWEEKS} that the Weeks manifold, with volume 0.9427..., is the hyperbolic 3--manifold of minimum volume. Equivalent results are unknown for higher dimensions. 

In this section, we reference several known results on volume, in terms of the hyperbolic metric, for several subcategories of hyperbolic $n$--orbifolds. For ease of comparison, we approximate to two significant digits. The bounds achieved in this paper improve upon the general hyperbolic manifold and orbifold bounds known to the authors \cite{Ade1}, \cite{FriHer}, \cite{Mar1}. However, they are smaller than the sharp bounds given for cusped and arithmetic orbifolds. 

\subsection{Hyperbolic Manifolds} 

An explicit lower bound for the radius of a ball that can be embedded in every hyperbolic $n$--manifold was given in \cite{Mar1}. An error in that paper was later corrected in \cite{FriHer}. Using the corrected radius, \[\frac{0.0025}{17^{\lfloor n/2\rfloor}},\] one can obtain a lower bound for the volume of a hyperbolic $n$--manifold. In dimension three the bound is $1.33\times 10^{-11}$. 

\subsection{Cusped Hyperbolic Orbifolds} 

The smallest cusped hyperbolic 3--orbifold has volume $7.22\times10^{-2}$ \cite{Mey2}. The bound is $6.85\times10^{-3}$ in dimension four \cite{HildKell}. Analogous results for all dimensions less than ten can be found in \cite{Hild}.

\subsection{Arithmetic Orbifolds} It is conjectured that, for each dimension, the hyperbolic orbifold (also manifold) of minimum volume is arithmetic. This is the case in dimensions two and three \cite{Chin1}, \cite{Chin2}. The minimal volume arithmetic $n$--orbifolds were identified for all dimensions greater than or equal to four in \cite{Belo}, \cite{BeloOops}, \cite{Belo2}. 

For example, the volume of the smallest compact arithmetic hyperbolic $n$--orbifold, for $n=2r$ and $r$ even, is given by
\[\omega_c(n)=\frac{4\cdot5^{r^2+r/2} \cdot (2\pi)^r}{(2r-1)!!} \prod_{i=1}^r \frac{(2i-1)!^2}{(2\pi)^{4i}}\zeta_{k_0}(2i).\]
Here, $\zeta_{k_0}$ represents the Dedekind zeta function of the number field $k_0=\mathbb Q[\sqrt 5]$. 

The cited papers contain similar formulas for $n=2r, r$ odd and $n=2r-1$. The noncompact cases are also addressed. The volume of the smallest compact arithmetic hyperbolic 4--orbifold, which is extremal among the volumes of all known hyperbolic 4--orbifolds, is calculated to be $1.8\times10^{-3}$.

 \section*{Acknowledgments} 
 
 This project was initiated during the first author's fellowship at the Mathematical Sciences Research Institute in the Fall of 2007. Most of this work was done while he held a visiting position at the University of California, Santa Barbara. In addition, he is grateful to Francis Bonahon, Dick Canary and Daryl Cooper for their support and useful conversations. 
 
 The second author is partially supported by NSF grant DMS-0806016.

\bigskip

\bibliography{ProjectTwo.bib}

\end{document}